\documentclass[12pt]{article}
\usepackage{mathptmx}
\usepackage{graphicx}
\usepackage{amsmath}
\usepackage{amsfonts}
\newtheorem{thm}{Theorem}[section]
\newtheorem{lem}{Lemma}[section]
\newtheorem{rmk}{Remark}[section]

\newtheorem{Corollary}{Corollary}[section]

\newtheorem{tb}{Table}
\usepackage{epstopdf}

\title{ Approximation by Durrmeyer type Exponential Sampling Series}
\author{Shivam Bajpeyi \thanks{Department of Mathematics, Visvesvaraya National Institute of Technology, Nagpur, Nagpur-440010, India. \newline E-mail: shivambajpai1010@gmail.com }
 \and
  A. Sathish Kumar \thanks{Department of Mathematics, Visvesvaraya National Institute of Technology, Nagpur, Nagpur-440010, India. \newline E-mail: mathsatish9@gmail.com}
  }
\date{}

\begin{document}
\maketitle
\bibliographystyle{plain}
\abstract{In this article, we analyze the approximation properties of the new family of Durrmeyer type exponential sampling operators. We derive the point-wise and uniform approximation theorem and Voronovskaya type theorem for these generalized family of operators. Further, we construct a convex type linear combination of these operators and establish the better approximation results. Finally, we provide few examples of the kernel functions to which the presented theory can be applied along with the graphical representation.

\endabstract



\section{Introduction}
The sampling theory provides the prominent tool to handle various problems arising in signal processing and approximation theory. The classical sampling theorem is credit to the names of Whittaker-Kotel'nikov-Shannon, which provides a reconstruction formula for the band-limited functions using its ordinates at a series of points (see \cite{wks,shan}). In order to weaken the assumption, Brown \cite{brown}, Butzer and Splettsst\"{o}ber \cite{spel} initiated the study and since then, many authors have contributed significantly in the development of the theory in various settings, see eg. \cite{butzer2,bardaro10,bardaro4,costa1,buta,vinti1,vinti2,bardaro5,ang}.
\par
The theory of exponential sampling began due to a group of engineers and physicists Pike, Bertero \cite{bertero} and Gori \cite{gori} who introduced the \textit{exponential sampling formula} (see also \cite{ostrowsky,butzer5}) which provides a series representation for the Mellin band limited functions exploiting its exponentially spaced sample values. The theory of exponential sampling has been proven to be an important tool to deal with various inverse problems in the area of optical physics, see eg. \cite{bertero,casasent,gori,ostrowsky}. The mathematical study of the exponential sampling theorem was firstly initiated by Butzer and Jansche in \cite{butzer5} using the Mellin theory. The Mellin theory revealed to be the most suitable frame to put the theory of exponential sampling and found many applications in the boundary value problems (see \cite{butzer5}). Mamedov pioneered the separate study of Mellin theory independent from Fourier analysis and investigated the approximation properties of Mellin convolution operators in \cite{mamedeo}. We mention here some of the significant developments in the Mellin analysis \cite{butzer3,butzer5,butzer7,bard,bardaro9,bardaro3}.
\par
Bardaro et.al.\cite{bardaro7} made a remarkable development in this direction by generalizing the exponential sampling formula, where they considered the generalized kernel with suitable assumptions in place of $lin_{c}$-function. This furnished a mechanism to approximate not necessarily Mellin band-limited functions using the sample values at nodes $(e^{\frac{k}{w}})_{w>0},\ k \in \mathbb{Z}.$ These operators have been studied in different settings in \cite{bardaro1,bardaro11,bardaro8,comboexp}. In order to reduce the \textit{time-jitter} error, the Kantorovich modification of these operators has been introduced and studied in \cite{own,lnr}. The \textit{time-jitter} error causes when the sample values can not be obtained exactly at the nodes. The Kantorovich type operators are known to reduce the \textit{time-jitter} error as they calculate the information around the nodes rather than exactly at the nodes. This modification enables to approximate integrable functions using the sample values which are exponentially spaced. The Kantorovich version of the several operators have been investigated extensively in various settings, see \cite{bardaro10,bardaro5,orlova,ang,vinti1,vinti2,sir,gbs,ana,costa1}. The Kantorovich type exponential sampling series was defined in \cite{own}, namely
\begin{equation} \label{kant}
(I_{w}^{\chi}f)(x)= \sum_{k= - \infty}^{+\infty} \chi(e^{-k} x^{w})\  w \int_{\frac{k}{w}}^{\frac{k+1}{w}} f(e^{u})\  du \ , \hspace{0.4cm} w>0
\end{equation}
where $ f: \mathbb{R}^{+} \rightarrow \mathbb{R}$ is locally integrable such that the above series is convergent for every $ x \in \mathbb{R}^{+}$.
\par
In view of Theorem 3.1 in \cite{lnr}, it is apparent that the Kantorovich exponential sampling operators (\ref{kant}) fails to improve the order of approximation in the asymptotic formula. This motivates us to adopt another approach, known as Durrmeyer method, where the integral mean is replaced by a general convolution operator. This method was firstly applied to the Bernstein polynomials in the series of papers \cite{D1,D3,D6}. The Durrmeyer type modification is known to provide the better order of approximation in various settings, see eg.\cite{ana,D7,D8,D9}.
In this paper, we generalize the operator (\ref{kant}) by replacing the integral mean with the Mellin singular integral to obtain the following family of operators. For $x \in \mathbb{R}^+$ and $w>0,$ we introduce the Durrmeyer type exponential sampling series as
\begin{equation*}
(I_{w}^{\chi,\phi}f)(x)= \sum_{k= - \infty}^{+\infty} \chi(e^{-k} x^w)\ w  \int_{0}^{\infty} \phi(e^{-k} t^{w}) f(t)\  \frac{dt}{t}.
\end{equation*}
where $ f: \mathbb{R}^{+} \rightarrow \mathbb{R}$ is integrable such that the above series is convergent for every $ x \in \mathbb{R}^{+}$.
\par
The paper is organized as follows. We begin with some auxiliary definitions in Section 3. In Section 4, we establish the basic convergence theorem and asymptotic formula for the proposed family operators (\ref{main}). Using Peetre's K functional \cite{anas,lnr}, we obtain the quantitative estimates of the Voronovskaya type theorem. Section 5 is devoted to the study of suitable linear combinations of the family of operators (\ref{main}) which produce the better order of approximation. The idea of considering the linear combinations of the operators is a prominent method to improve the order of convergence, see \cite{lc1,lc2,comboexp,D7,D8,lnr}. In last section, we discuss few examples of the well known kernels along with the graphical representations.

\section{Preliminaries}
Let $\mathbb{R}^+$ be the multiplicative topological group equipped with the Haar measure $\displaystyle \mu(S)= \int_{S} \frac{dt}{t},$ where $dt$ represents the Lebesgue measure and $S$ is any measurable set. We define $L^{p}(\mu,\mathbb{R})=: L^{p}(\mu),\ 1 \leq p < +\infty,$ the Lebesgue spaces with respect to the measure $\mu,$ endowed with the usual $p$- norm. We consider $M(\mu)$ as the class of all measurable functions  and $L^{\infty}(\mu)$ as the space of all bounded functions with respect to the measure $\mu.$
\par
Let $X_{c}$ denotes the space of all functions $f: \mathbb{R}^+ \rightarrow \mathbb{R}$ such that $f(x) x^{c} \in L^{1}(\mu)$ for some $c \in \mathbb{R},$ endowed with the norm
$$\Vert f \Vert_{X_c} = \int_0^{+\infty} |f(u)|u^{c}\  \frac{du}{u}.$$
The Mellin transform of $f \in X_{c}$ is defined by
$$\hat{M}[f](s) := \int_0^{+\infty} f(u) u^{s}\ \frac{du}{u} , \ \ \ (s = c + it, \ t \in \mathbb{R}).$$ Indeed, if the Mellin transform exists for some $c$ then it exists for every $s=c+it.$ The basic properties of the Mellin transform can be found in \cite{butzer3,bardaro7}.
Furthermore, the pointwise derivative in the Mellin's frame is given by the following limit:
$$ \lim_{h \rightarrow 1} \frac{\tau_{h}^{c}f(x) - f(x)}{h-1} = x f^{'}(x)+c f(c),$$ provided $f^{'}$ exists. Here, $\tau_{h}^{c}$ represents the Mellin translation operator and is defined as $(\tau_{h}^{c}f)(x):= h^{c} f(hx).$ Thus the pointwise Mellin's derivative $\theta_{c}f$ of any function $f : \mathbb{R}^+ \rightarrow \mathbb{C}$ defined as
$$(\theta_cf)(x) := xf'(x) + cf(x), \ \ \ \ \ x \in  \mathbb{R}^+$$ provided $f^{'}$ exists.
Subsequently, the $r^{th}$ order Mellin differential operator can be expressed as $$\theta_c^r:=\theta_c(\theta_c^{r-1}).$$ For the sake of convenience, we define $\theta_c:=\theta_c^1$ and $\theta f:=\theta_{0} f.$
\vspace{0.2cm}

Moreover, $C(\mathbb{R}^+)$ denotes the space of all continuous and bounded functions on $\mathbb{R}^+$ equipped with the supremum norm $\|f\|_{\infty} := \sup_{x \in \mathbb{R}^+} |f(x)|.$ Subsequently, for any $r \in \mathbb{N},$ $C^{(r)}(\mathbb{R}^+)$ be the subspace of $C(\mathbb{R}^+)$ such that $f^{(k)},\ k \in \mathbb{N}$ exists for every $k \leq r$ and each $f^{(k)} \in C(\mathbb{R}^+).$ In what follows, we call a function $f: \mathbb{R}^+ \rightarrow \mathbb{C}$ log-uniformly continuous on $\mathbb{R}^+$ if for any given $\epsilon > 0,$ there exists $\delta > 0$ such that $|f(u) -f(v)| < \epsilon$ whenever $| \log u - \log v| < \delta,$ for any $u, v \in \mathbb{R}^{+}.$ We denote $\mathcal{C} (\mathbb{R}^+)$ the space of $C(\mathbb{R}^+)$ containing all log-uniformly continuous and bounded functions defined on  $\mathbb{R}^{+}.$ Similarly, $\mathcal{C}^{(r)}(\mathbb{R}^+)$ denotes the space of functions which are $r$-times continuously Mellin differentiable and $\theta^{r}f \in \mathcal{C}(\mathbb{R}^+).$ The notion of log-continuity was first introduced in \cite{mamedeo}.

\section{Durrmeyer type generalization of exponential sampling series}
Let $ \chi :\mathbb{R}^{+} \rightarrow \mathbb{R}$ be the kernel which is continuous on $\mathbb{R}^{+}.$ For any $\nu \in \mathbb{N}_{0}:= \mathbb{N} \cup \{0\}$ and $x \in \mathbb{R}^{+},$ we define its algebraic and absolute moments of order $\nu$ respectively as
$$ m_{\nu}(\chi,u):= \sum_{k= - \infty}^{+\infty}  \chi(e^{-k} u) (k- \log u)^{\nu}, \hspace{0.5cm} M_{\nu}(\chi):= \sup_{u \in \mathbb{R}^+} \sum_{k= - \infty}^{+\infty}  |\chi(e^{-k} u)| |k- \log u|^{\nu}.$$
Moreover, let $\phi$ be the kernel such that $\phi \in L^{1}(\mu).$ Then the algebraic and absolute moments for the kernel $\phi$ are defined by
$$ \hat{m}_{\nu}(\phi):= \int_{0}^{\infty} \phi(u) \log^{\nu} u \ \frac{du}{u}, \hspace{0.4cm} \hat{M}_{\nu}(\phi):= \int_{0}^{\infty} |\phi(u)| \  |\log u|^{\nu} \ \frac{du}{u}.$$
We suppose that the kernel functions $ \chi$ and $\phi$ satisfy the following conditions:
\begin{itemize}
\item[K1)] For every $ \displaystyle u \in \mathbb{R}^{+}, \hspace{0.3cm}\sum_{k=- \infty}^{+\infty} \chi(e^{-k} u) =1 \ \ \ \ \mbox{and} \hspace{0.5cm} \int_{0}^{\infty} \phi(u) \frac{du}{u}=1.$
\item[K2)] For some $r \in \mathbb{N},$ $(M_{r}(\chi) +\hat{M}_{r}(\phi)) <+ \infty$ and
$$\displaystyle \lim_{\gamma \rightarrow + \infty} \sum_{|k-\log u|>  \gamma} |\chi(e^{-k} u)| \ |k- \log u|^{r}=0$$ uniformly with respect to $ u \in \mathbb{R}^{+}.$
\end{itemize}

\begin{rmk} \cite{own}
It is easy that for $\mu, \nu \in \mathbb{N}_{0}$ with $\mu < \nu,$ $M_{\nu}(\chi) < + \infty$ implies that $M_{\mu}(\chi) < + \infty.$ Moreover, the condition (K2) implies that there holds
$$\lim_{\gamma \rightarrow + \infty} \sum_{|k-\log u|>  \gamma} |\chi(e^{-k} u)| \ |k- \log u|^{j}=0 \ \ \mbox{for}\ \ j=0,1,...,r-1 \ .$$
\end{rmk}
Following along the lines of Remark 4.1 in \cite{own}, we deduce that $\hat{M}_{\nu}(\phi) < + \infty$ implies that $\hat{M}_{\mu}(\phi) < + \infty$ whenever $\mu < \nu.$

\vspace{0.3cm}
Under the above assumptions on the kernels $\chi$ and $\phi,$ we define the Durrmeyer type generalized exponential sampling series as
\begin{equation} \label{main}
(I_{w}^{\chi,\phi}f)(x)= \sum_{k= - \infty}^{+\infty} \chi(e^{-k} x^w)\ w  \int_{0}^{\infty} \phi(e^{-k} t^{w}) f(t)\  \frac{dt}{t}\ , \hspace{0.4cm} w>0
\end{equation}
where $ f: \mathbb{R}^{+} \rightarrow \mathbb{R}$ is integrable such that the above series is convergent for every $ x \in \mathbb{R}^{+}$. It is evident that for $f \in L^{\infty}(\mu),$ the series (\ref{main}) is well-defined, that is, $ L^{\infty} (\mu) \in Dom (I_{w}^{\chi,\phi}),$ where $Dom (I_{w}^{\chi,\phi})$ consists of all functions $f: \mathbb{R}^{+} \rightarrow \mathbb{R}$ such that the series (\ref{main}) is absolutely convergent for any $x \in \mathbb{R}^+.$

\begin{rmk}
Let $\chi$ and $\phi$ be the kernel functions. Then the Mellin convolution integral of $f$ with $\phi_{w}(u) := w\  \phi(u^{w})$ is given by
$$(T_{w}^{\phi}f)(s):= \big(f*\phi_{w}(u) \big)(s)= \int_{0}^{\infty}f(t)\  \phi_{w} \left( \frac{t}{s} \right) \frac{dt}{t}=\  w\int_{0}^{\infty}f(t)\  \phi\left( \frac{t^w}{s^w} \right).$$
Now we can represent our proposed operator (\ref{main}) in terms of the generalized exponential sampling operator $(S_{w}^{\chi})_{w>0}$ introduced by Bardaro et.al. in \cite{bardaro7}, using the Mellin convolution integral as follows
$$(I_{w}^{\chi,\phi}f)(x):= (S_{w}^{\chi}\left(T_{w}^{\phi}f \right))(x)\ , \hspace{0.5cm} x \in \mathbb{R}^+. $$
\end{rmk}

\begin{rmk}
For any $y \in \mathbb{R}^+,$ we define
\begin{equation*}
\phi(y) := \kappa_{[1, e]}(y) =
     \begin{cases}

     {1 ,} &\quad\text{} \ \ \ \ { 1 \leq y < e}\\
       {0,} &\quad\text{} \ \ \ \ { otherwise }\\
   \end{cases}
\end{equation*}
where $\kappa$ represents the characteristic function. Considering the change of variable $t=e^u,$ we obtain
\begin{equation*}
\phi(e^{wu-k}) =
     \begin{cases}

     {1 ,} &\quad\text{} \ \ \ \ { \frac{k}{w} \leq u \leq \frac{k+1}{w} }\\
       {0,} &\quad\text{} \ \ \ \ { otherwise }.\\
   \end{cases}
\end{equation*}
Then, the corresponding generalized operator (\ref{main}) acquires the following form:
\begin{eqnarray} \label{expo}
(I_{w}^{\chi}f)(x)= \sum_{k= - \infty}^{+\infty} \chi(e^{-k} x^w)\ w  \int_{\frac{k}{w}}^{\frac{k+1}{w}} f(e^{u}) \ du.
\end{eqnarray}
This particular family of operators (\ref{expo}) was analyzed in \cite{own,lnr}.
\end{rmk}

\section{Convergence theorems}
The purpose of this section is to derive some local approximation results for the Durrmeyer type exponential sampling operators (\ref{main}).

\begin{thm}\label{theorem1}
Let $ f\in M(\mu) \cap L^{\infty}(\mu).$ Then the series (\ref{main}) converges to $f(x)$ at every point $x \in \mathbb{R}^{+},$ the point of continuity of $f$. Moreover, for $ f \in C(\mathbb{R}^{+}),$ we have
$$ \lim_{w \rightarrow \infty} \|I_{w}^{\chi,\phi}f - f \|_{\infty} = 0.$$
\end{thm}

\noindent\bf{Proof.}\rm \
Using the condition (K1), we can write
\begin{eqnarray*}
|(I_{w}^{\chi,\phi}f)(x)-f(x)| &=& \bigg| \sum_{k= - \infty}^{+\infty} \chi(e^{-k} x^w)\ w \int_{0}^{\infty}  \phi(e^{-k} t^{w}) (f(t) - f(x)) \frac{dt}{t} \bigg| \\
&\leq & \sum_{k= - \infty}^{+\infty} \big| \chi(e^{-k} x^w) \big| w \left( \int_{|\log t - \log x |< \delta} +\int_{|\log t - \log x | \geq \delta} \right) \big| \phi(e^{-k} t^{w})\big|\ \  \big| f(t) - f(x) \big| \frac{dt}{t} \\
&:=& I_{1}+I_{2}. \
\end{eqnarray*}
Let $x \in \mathbb{R}^+$ be the point of continuity of $f$ then for any fixed $ \epsilon >0$ there exists $\delta>0$ such that $ |f(t) - f(x)|< \epsilon,$ whenever $|\log t - \log x|< \delta.$ Now, considering the change of variable $e^{-k} t^{w}=:p, $ we obtain
\begin{eqnarray*}
|I_{1}| & \leq & \sum_{k= - \infty}^{+\infty} \big| \chi(e^{-k} x^w) \big| \ \epsilon \ \int_{w (\log x-\delta) -k} ^{w ( \log x+\delta) -k} \big| \phi(p)\big| \ \frac{dp}{p} \\
& \leq & \epsilon  \big( M_{0}(\chi) \ \|\phi\|_{1} \big).
\end{eqnarray*}
Similarly, we estimate $ I_{2}.$
\begin{eqnarray*}
|I_{2}| & \leq &  2 \|f\|_{\infty} \sum_{k \in \mathbb{Z}} \big| \chi(e^{-k} x^w)\big| \left( \int_{-\infty}^{w(\log x - \delta)-k} + \int_{w(\log x + \delta)-k}^{+\infty} \right) |\phi(p)| \frac{dp}{p}\\
&=:& I_{2}^{'}+I_{2}^{''}.
\end{eqnarray*}
First we consider $I_{2}^{'}.$ Using the fact that $\phi \in L^{1}(\mu)$ and condition (K2), for a fixed $x, \delta$ we obtain
\begin{eqnarray*}
|I_{2}^{'}|  & \leq & 2 \|f\|_{\infty}  \left( \sum_{|w \log x - k|< \frac{w \delta}{2}} +\sum_{|w \log x-k| \geq \frac{w \delta}{2}} \right) \big| \chi(e^{-k} x^w)\big| \int_{-\infty}^{w(\log x - \delta)-k} |\phi(p)| \frac{dp}{p} \\
& \leq & 2 \|f\|_{\infty} \epsilon (1+ \|\phi\|_{1}).
\end{eqnarray*}
Similarly, we deduce that $|I_{2}^{''}| \leq 2 \|f\|_{\infty} \epsilon (1+ \|\phi\|_{1}).$ Finally combining the estimates $I_{1}-I_{2},$ we get the desired result. Subsequently for $f \in C(\mathbb{R}^+),$ the proof follows in the similar manner.
\vspace{0.2cm}

Our next result is the following asymptotic formula for the family of operators $(I_{w}^{\chi,\phi})_{w>0} .$
\begin{thm}
Let $ f \in C^{(r)}(\mathbb{R}^{+})$ locally at $x \in \mathbb{R}^+$ and  $\chi, \phi$ be the kernels. Then we have
$$ \big[(I_{w}^{\chi,\phi}f)(x) - f(x) \big] = \sum_{j=1}^{r} \frac{\theta^{j}f(x)}{j! \ w^{j}} \left( \sum_{\eta=0}^{j} \binom{j}{\eta} \hat{m}_{j-\eta}(\phi)\  m_{\eta}(\chi,x) \right)+R_{w,r}(x),$$
where $\displaystyle R_{w,r}(x)= \sum_{k \in \mathbb{Z}} \chi(e^{-k} x^w) \ \  w \int_{0}^{\infty} \phi(e^{-k} t^{w}) \ h \left( \frac{t}{x} \right) (\log t - \log x)^{r} \ \frac{dt}{t}.$ \\
Moreover $R_{w,r}(x)= o(w^{-r})$ as $w \rightarrow +\infty$ for every $x \in \mathbb{R}^+.$
\end{thm}

\noindent\bf{Proof.}\rm \ Using Taylor's formula in terms of Millan derivatives (see \cite{bardaro7,own}), we write
$$\displaystyle f(t)  - f(x)= \sum_{j=1}^{r} \frac{\theta^{j}f(x)}{j!}(\log t - \log x)^{j} + h \left( \frac{t}{x} \right) (\log t - \log x)^{r},$$ where $h$ is a bounded function such that $\displaystyle \lim_{y \rightarrow 1} h(y)=0.$ In view of (\ref{main}) we obtain
\noindent \\
$[(I_{w}^{\chi,\phi}f)(x) - f(x)]$
\begin{eqnarray*}
&=& \sum_{k= - \infty}^{+\infty} \chi(e^{-k} x^w)\ w \int_{0}^{\infty}  \phi(e^{-k} t^{w}) \Bigg[ \sum_{j=1}^{r} \frac{\theta^{j}f(x)}{j!}(\log t - \log x)^{j} + h \left( \frac{t}{x} \right) (\log t - \log x)^{r} \Bigg] \ \frac{dt}{t} \\
&=& \sum_{j=1}^{r} \frac{\theta^{j}f(x)}{j!} \sum_{k= - \infty}^{+\infty} \chi(e^{-k} x^w)\ w \int_{0}^{\infty}  \phi(e^{-k} t^{w})(\log t - \log x)^{j} +R_{w,r}(x),
\end{eqnarray*}
where $\displaystyle R_{w,r}(x)= \sum_{k= - \infty}^{+\infty} \chi(e^{-k} x^w)\ w \int_{0}^{\infty}  \phi(e^{-k} t^{w}) \ h \left( \frac{t}{x} \right) (\log t - \log x)^{r} \ \frac{dt}{t}.$
\vspace{0.3cm}
\\
For any fixed index $\xi \in \mathbb{N},$ we have
\noindent \\

$\displaystyle \sum_{k= - \infty}^{+\infty} \chi(e^{-k} x^w)\ w \int_{0}^{\infty}  \phi(e^{-k} t^{w}) \frac{\theta^{\xi}f(x)}{\xi !}(\log t - \log x)^{\xi} \ \frac{dt}{t}$
\begin{eqnarray*}
 &=& \frac{\theta^{\xi}f(x)}{\xi !} \sum_{k= - \infty}^{+\infty} \chi(e^{-k} x^w)\ w \int_{0}^{\infty}  \phi(e^{-k} t^{w})(\log t - \log x)^{\xi} \ \frac{dt}{t} \\
&=& \frac{\theta^{\xi}f(x)}{\xi! \ w^{\xi}} \sum_{k= - \infty}^{+\infty} \chi(e^{-k} x^w) \int_{0}^{\infty} \phi(p)(\log p+ k - w \log x)^{\xi} \ \frac{dp}{p}\\
&=& \frac{\theta^{\xi}f(x)}{\xi! \ w^{\xi}} \sum_{k= - \infty}^{+\infty} \chi(e^{-k} x^w) \int_{0}^{\infty} \phi(p) \left( \sum_{\eta=0}^{\xi} \binom{\xi}{\eta} (\log p)^{(\xi-\eta)} (k-w \log x)^{\eta} \right) \ \frac{dp}{p}\\
&=& \frac{\theta^{\xi}f(x)}{\xi! \ w^{\xi}} \left( \sum_{\eta=0}^{\xi} \binom{\xi}{\eta} \hat{m}_{\xi-\eta}(\phi) \ m_{\eta}(\chi,x) \right).
\end{eqnarray*}
Thus we obtain
\begin{eqnarray} \label{eq}
[(I_{w}^{\chi,\phi}f)(x) - f(x)] &=& \sum_{j=1}^{r} \frac{\theta^{j}f(x)}{j! \ w^{j}} \left( \sum_{\eta=0}^{j} \binom{j}{\eta} \hat{m}_{j-\eta}(\phi)\  m_{\eta}(\chi,x) \right)+ R_{w,r}(x).
\end{eqnarray}
Now we estimate the remainder term $R_{w,r}(x)$ in (\ref{eq}). In order to do that let $\epsilon > 0$ be fixed then there exists $\delta > 0$ such that $|h(y)| < \epsilon $ whenever $|\log y|< \delta.$ We write $R_{w,r}(x)$ as
\begin{eqnarray*}
|R_{w,r}| &\leq & \sum_{k \in \mathbb{Z}} \big |\chi(e^{-k} x^w)\big | \ \ \Bigg| w \int_{|\log t- \log x| < \delta} \phi(e^{-k} t^{w}) \ h \left( \frac{t}{x} \right) (\log t - \log x)^{r} \ \frac{dt}{t} \Bigg| \\&& + \sum_{k \in \mathbb{Z}} \big|\chi(e^{-k} x^w)\big| \ \ \Bigg| w  \int_{|\log t - \log x| \geq \delta} \phi(e^{-k} t^{w}) \ h \left( \frac{t}{x} \right) (\log t - \log x)^{r} \ \frac{dt}{t} \Bigg| \\
&:=& I^{'}+I^{''}.
\end{eqnarray*}
Thus we have
$$|I^{'}| \leq \frac{\epsilon}{w^{r}} \left ( \sum_{j=0}^{r} \binom{r}{j} M_{(r-j)}(\chi) \hat{M}_{j}(\phi) \right).$$
Using the fact that $h$ is bounded, we obtain
\begin{eqnarray*}
|I^{''}| &\leq & \|h\|_{\infty} \sum_{k \in \mathbb{Z}} \big |\chi(e^{-k} x^w)\big | \ \ \Bigg| w \int_{|\log t - \log x| \geq \delta} \phi(e^{-k} t^{w}) (\log t - \log x)^{r} \frac{dt}{t} \Bigg| \\
&\leq & \frac{\|h\|_{\infty}}{w^r}  \sum_{k \in \mathbb{Z}}  \big |\chi(e^{-k} x^w)\big | \int_{|\log p+k-w \log x| \geq w \delta} | \phi(p)| |\log p +k-w \log x|^{r} \frac{dp}{p} \\
&\leq & \frac{2^{r-1} \epsilon \|h\|_{\infty}}{w^r} \left(M_{r}(\chi) + \| \phi \|_{1} \right).
\end{eqnarray*}
On Combining the estimates of $I^{'}-I^{''},$ we conclude that $R_{w,r}(x)= o(w^r)$ as $w \rightarrow +\infty$ and hence the desired result is established.\\

From the above theorem, we deduce the following Voronovskaja type asymptotic result.
\begin{Corollary}
Let $ f \in C^{(1)}(\mathbb{R}^{+})$ and  $\chi, \phi$ be the kernels. Then for every $x \in \mathbb{R}^+,$
$$ \lim_{w \rightarrow \infty} w \big[(I_{w}^{\chi,\phi}f)(x) - f(x) \big] = \theta f(x) \left( m_{1}(\chi,x) +\hat{m}_{1}(\phi) \right).$$
\end{Corollary}

\begin{Corollary}
Under the assumptions of Theorem 4.2, if in addition, the kernels $\chi, \phi$ satisfy $m_{j}(\chi,x)=0$ and $\hat{m}_{j}(\phi)=0$ for $j=1,2,...\ r-1,$ then the following holds
$$ \lim_{w \rightarrow \infty} w^{r} [(I_{w}^{\chi,\phi}f)(x) - f(x)]= \frac{\theta^{r}f(x)}{r!} (m_{r}(\chi,x)+ \hat{m}_{r}(\phi)).$$
\end{Corollary}

\begin{rmk}
The condition that f is bounded, can be relaxed by assuming that $|f(e^x)| \leq(a+ b|x|), \ x \in \mathbb{R}^+,$ where $a,b \in \mathbb{R}$ are arbitrary constants.
Indeed, from (\ref{main}), we have
\begin{eqnarray*}
|(I_{w}^{\chi,\phi}f)(x)| & \leq & \sum_{k= - \infty}^{+\infty} |\chi(e^{-k} x^w)| \  w \int_{0}^{\infty} |\phi (e^{-k} t^{w}) | \  (a+b |\log t|) \frac{dt}{t} \\
& \leq & \sum_{k= - \infty}^{+\infty} |\chi(e^{-k} x^w)| \int_{0}^{\infty} |\phi(p)| \left( a+ \frac{b}{w} \big |\log p +k \big | \right) \frac{dp}{p}\\
& \leq & (a+ b |\log x|) M_{0}(\chi) \hat{M}_{0}(\phi) + \frac{b}{w} \left( M_{0}(\chi) \hat{M}_{1}(\phi) + M_{1}(\chi) \hat{M}_{0}(\phi)\right).
\end{eqnarray*}
Under the assumptions on the kernels $\chi$ and $\phi,$ we conclude that $f \in Dom(I_{w}^{\chi,\phi}).$
\end{rmk}
\subsection{Quantitative estimates}
We establish the quantitative estimate of the asymptotic formula derived in Theorem 4.2, in terms of Peetre's K-functional \cite{anas,lnr}. The Peetre's K-functional for $f \in \mathcal{C}^{(r)}(\mathbb{R}^+)$ is defined by
$$ \hat{K}(f,\epsilon, \mathcal{C}^{(r)}(\mathbb{R}^+),\mathcal{C}^{(r+1)}(\mathbb{R}^+)) := \inf \{ \|\theta^{r}(f-g)\|_{\infty} + \epsilon \|\theta^{r+1} g \|_{\infty} : g \in \mathcal{C}^{(1)}(\mathbb{R}^+), \epsilon \geq 0 \} .$$

\begin{thm}
Let $\chi, \phi$ be the kernel functions and $f \in \mathcal{C}^{(r)}(\mathbb{R}^+).$ Then we have
$$ \Bigg| [(I_{w}^{\chi}f)(x) - f(x)] - \sum_{j=1}^{r} \frac{\theta^{j}f(x)}{j! \ w^{j}} \left( \sum_{\eta=0}^{j} \binom{j}{\eta} \hat{m}_{j-\eta}(\phi)\  m_{\eta}(\chi,x) \right) \Bigg| \leq \frac{2^r \ A}{r! w^r}\hat{K} \left( \theta^{r}f, \frac{1}{(r+1)w} \frac{B}{A} \right).$$
where $A:= \left( M_{0}(\chi) \hat{M}_{r}(\phi)+M_{r}(\chi) \hat{M}_{0}(\phi) \right)$ and $B:=\left( M_{0}(\chi) \hat{M}_{r+1}(\phi)+ M_{r+1}(\chi) \hat{M}_{0}(\phi) \right).$
\end{thm}

\noindent\bf{Proof.}\rm \ From Theorem 3.2, we have
\\

\noindent $\displaystyle \Bigg| [(I_{w}^{\chi}f)(x) - f(x)] - \sum_{j=1}^{r} \frac{\theta^{j}f(x)}{j! \ w^{j}} \left( \sum_{\eta=0}^{j} \binom{j}{\eta} \hat{m}_{j-\eta}(\phi)\  m_{\eta}(\chi,x) \right) \Bigg|$
\begin{eqnarray*}
 =\displaystyle \Bigg| \sum_{k= - \infty}^{+\infty} \chi(e^{-k} x^w)\ w \int_{0}^{\infty}  \phi(e^{-k} t^{w}) h \left( \frac{t}{x} \right) (\log t - \log x)^{r} \frac{dt}{t}\Bigg|.
\end{eqnarray*}

Now we substitute $\displaystyle R_{r}(f,x,t):= h \left( \frac{t}{x} \right) (\log t - \log x)^{r}.$ Then by using the following estimate (see \cite{bardaro9,lnr}) $$ \displaystyle |R_{r}(f,x,t)| \leq \frac{2}{r!} |\log t - \log x|^{r} \ \hat{K} \left( \theta^{r}f, \frac{|\log t- \log x|}{2(r+1)} \right),$$  we can write \\

\noindent $\displaystyle \Bigg| [(I_{w}^{\chi}f)(x) - f(x)] - \sum_{j=1}^{r} \frac{\theta^{j}f(x)}{j! \ w^{j}} \left( \sum_{\eta=0}^{j} \binom{j}{\eta} \hat{m}_{j-\eta}(\phi)\  m_{\eta}(\chi,x) \right) \Bigg|$
\begin{eqnarray*}
&\leq & \sum_{k= - \infty}^{+\infty} |\chi(e^{-k} x^w)| \ w \int_{0}^{\infty}  |\phi(e^{-k} t^{w})| \ \left(\frac{2}{r!} |\log t - \log x|^r \right) \ \hat{K} \left( \theta^{r}f, \frac{|\log t - \log x|}{2(r+1)} \right) \frac{dt}{t}\\
 & \leq & \sum_{k= - \infty}^{+\infty} |\chi(e^{-k} x^w)| \ (2w) \int_{0}^{\infty}  |\phi(e^{-k} t^{w})| \ \frac{|\log t - \log x|^r}{r!} \left( \|\theta^{r}(f-g)\|_{\infty} + \frac{|\log t - \log x|}{2(r+1)} \|\theta^{r+1} g \|_{\infty} \right)\frac{dt}{t}\\
&\leq & \frac{2 \|\theta^{r}(f-g)\|_{\infty}}{r!} \sum_{k= - \infty}^{+\infty} |\chi(e^{-k} x^w)|  \ w \int_{0}^{\infty}  |\phi(e^{-k} t^{w})| \ |\log t - \log x|^r \frac{dt}{t} \ + \\&&
\frac{\|\theta^{r+1} g \|_{\infty}}{(r+1)!} \sum_{k= - \infty}^{+\infty} |\chi(e^{-k} x^w)|  \ w \int_{0}^{\infty}  |\phi(e^{-k} t^{w})| \ |\log t - \log x|^{r+1} \frac{dt}{t}\\
&:= & I_{1}+I_{2}.
\end{eqnarray*}
First we evaluate $I_{1}.$
\begin{eqnarray*}
|I_{1}| & \leq & \frac{2 \|\theta^{r}(f-g)\|_{\infty}}{r! \ w^{r} } \sum_{k= - \infty}^{+\infty} |\chi(e^{-k} x^w)| \int_{0}^{\infty}  |\phi(p)| |\log p+k-w \log x|^r \ \frac{dp}{p} \\
& \leq & \frac{2^r \|\theta^{r}(f-g)\|_{\infty}}{r! \ w^{r} } \left( \sum_{k= - \infty}^{+\infty} |\chi(e^{-k} x^w)| \int_{0}^{\infty}  |\phi(p)| |\log p|^r \ \frac{dp}{p} \right) + \\&&
 \frac{2^r \|\theta^{r}(f-g)\|_{\infty}}{r! \ w^{r} } \left( \sum_{k= - \infty}^{+\infty} |\chi(e^{-k} u)| |k - w\log x|^r \int_{0}^{\infty}  |\phi(p)| \ \frac{dp}{p} \right) \\
& \leq & \frac{2^r \|\theta^{r}(f-g)\|_{\infty}}{r! \ w^{r} } \left( M_{0}(\chi) \hat{M}_{r}(\phi)+M_{r}(\chi) \hat{M}_{0}(\phi) \right).
\end{eqnarray*}
Similarly, we have
\begin{eqnarray*}
|I_{2}| & \leq & \frac{\|\theta^{r+1} g \|_{\infty}}{(r+1)!\  w^{r+1}} \sum_{k= - \infty}^{+\infty} |\chi(e^{-k} x^w)|  \int_{0}^{\infty}  |\phi(p)|\ |\log p +k-w \log x|^{r+1} \frac{dp}{p}\\
& \leq & \frac{2^{r}\|\theta^{r+1} g \|_{\infty}}{(r+1)!\ w^{r+1}} \left( M_{0}(\chi) \hat{M}_{r+1}(\phi)+  M_{r+1}(\chi) \hat{M}_{0}(\phi)\right).
\end{eqnarray*}
On combining the estimates $I_{1}- I_{2}$ and taking the infimum to $g \in \mathcal{C}^{(r+1)}(\mathbb{R}^+),$ we establish the proof.
As a consequence of the above result, we deduce the following corollary.
\begin{Corollary}
Let $\chi, \phi$ be the kernel functions and $f \in \mathcal{C}^{(1)}(\mathbb{R}^+).$ Then the following estimate holds
$$ \Big| [(I_{w}^{\chi}f)(x) - f(x)] - (\theta f)(x) \big(\hat{m}_{1}(\phi)+  m_{1}(\chi) \big) \Big| \leq
\frac{2A}{w}\hat{K} \left( f, \frac{1}{2w} \frac{B}{A} \right),$$
where $A:= \left( M_{0}(\chi) \hat{M}_{1}(\phi)+M_{1}(\chi) \hat{M}_{0}(\phi) \right)$ and $B:=\left( M_{0}(\chi) \hat{M}_{2}(\phi)+ M_{2}(\chi) \hat{M}_{0}(\phi) \right).$
\end{Corollary}

\section{Improved order of approximation}
The aim of this section is to formulate the combinations of the family of operators (\ref{main}) to produce better order of convergence in the asymptotic formula for $f \in C^{(r)}(\mathbb{R}^+)$ without using the constraint that the higher order algebraic moments for the kernels $\chi, \phi$ vanish on $\mathbb{R}^+.$ We construct the linear combination of the family $(I_{w}^{\chi,\phi}(f,.))_{w >0}$ in the following manner.
\par
Let $\beta_{i}, \ i=1,2,...,p$ be non-zero real numbers such that $\displaystyle \sum_{i=1}^{p} \beta_{i}=1.$ For $x \in \mathbb{R}^+$ and $w >0,$ we define the linear combination of the operators (\ref{main}) as
\begin{eqnarray} \label{combo} \nonumber
(I_{w,p}^{\chi,\phi})(f,x)&=&\sum_{i=1}^{p} \beta_{i}  \sum_{k= - \infty}^{+\infty} \chi(e^{-k} x^{iw})\  (iw)  \int_{0}^{\infty}\phi(e^{-k} t^{iw}) f(t)\  \frac{dt}{t}\\
&=& \sum_{i=1}^{p} \beta_{i} \ (I_{iw}^{\chi,\phi})(f,x).
\end{eqnarray}
Now we derive the asymptotic formula for the family of operators defined in (\ref{combo}).

\begin{thm}\label{t4}
Let $\chi,\phi$ be the kernels and $f \in C^{(r)}(\mathbb{R}^+)$ locally at $x \in \mathbb{R}^+.$ Then we have
$$ [(I_{w,p}^{\chi,\phi})(f,x)- f(x)]= \sum_{j=1}^{r} \frac{(\theta^{j} f)(x)}{j!\ w^{j}} M_{j}^{p}(\chi,\phi)+o(w^{-r}),$$
where $\displaystyle M_{j}^{p}(\chi,\phi):= \sum_{i=1}^{p} \frac{\beta_{i}}{i^j} \left( \sum_{\eta=0}^{j} \binom{j}{\eta} \hat{m}_{j-\eta}(\phi) m_{\eta}(\chi,x)\right).$
\end{thm}
\noindent\bf{Proof.}\rm \
From the condition $ \displaystyle\sum_{i=1}^{p} \beta_{i}=1,$ we can write
$$ [(I_{w,p}^{\chi})(f,x)- f(x)]= \sum_{i=1}^{p} \beta_{i} \sum_{k= - \infty}^{+\infty} \chi(e^{-k} x^{iw})\  (iw)  \int_{0}^{\infty} \phi(e^{-k} t^{w}) (f(t) - f(x))\  \frac{dt}{t}. $$
Now using the $r^{th}$ order Mellin's Taylor formula, we obtain \\

\noindent $[(I_{w,p}^{\chi})(f,x)- f(x)]$
\begin{eqnarray*}
&=& \sum_{i=1}^{p} \beta_{i} \sum_{k= - \infty}^{+\infty} \chi(e^{-k} x^{iw})\ (iw)  \int_{0}^{\infty} \phi(e^{-k} t^{iw}) \left( \sum_{j=1}^{r}\frac{(\theta^{j} f)(x)}{j !} (\log t- \log x)^{j} + h \left( \frac{t}{x} \right) (\log t - \log x)^r \right)\ \frac{dt}{t} \\
&=& \sum_{i=1}^{p} \beta_{i} \sum_{k= - \infty}^{+\infty} \chi(e^{-k} x^{iw})\  (iw)  \int_{0}^{\infty} \phi(e^{-k} t^{iw}) \left( \sum_{j=1}^{r}\frac{(\theta^{j} f)(x)}{j !} (\log t- \log x)^{j} \frac{dt}{t} \right)\\&& +
 \sum_{i=1}^{p} \beta_{i} \sum_{k= - \infty}^{+\infty} \chi(e^{-k} x^{iw})\  (iw)  \int_{0}^{\infty} \phi(e^{-k} t^{iw}) \ h \left( \frac{t}{x} \right) (\log t- \log x)^r \ \frac{dt}{t} \\
&=& \sum_{i=1}^{p} \beta_{i} \left( \sum_{j=1}^{r}\frac{(\theta^{j} f)(x)}{j!\ w^{j} \ i^{j}} \left( \sum_{\eta=o}^{j} \binom{j}{\eta} (\hat{m}_{j-\eta}(\chi,x) m_{j}(\phi) \right) \right)+ R_{iw,r},
\end{eqnarray*}
where $\displaystyle R_{iw,r}= \sum_{i=1}^{p} \beta_{i} \sum_{k= - \infty}^{+\infty}\chi(e^{-k} x^{iw})\  (iw)  \int_{0}^{\infty} \phi(e^{-k} t^{iw}) h \left( \frac{t}{x} \right) (\log t- \log x)^r \ \frac{dt}{t}.$ Now proceeding along the lines of Theorem 4.2, it follows that $R_{iw,r}= o (w^{-r})$ as $w \rightarrow \infty.$ Hence, the proof is completed. \\

From the above result, we establish the following Voronovskaja type asymptotic formula.
\begin{Corollary}
For $f \in C^{(1)}(\mathbb{R}^+),$ we have
$$ w[(I_{w,p}^{\chi,\phi})(f,x)- f(x)] = \sum_{i=1}^{p} \frac{\beta_{i}}{i} \left(\theta f(x) (\hat{m}_{1}(\phi)+ m_{1}(\chi,x)) \right)+ o (w^{-1}).$$ Furthermore, if $m_{j}(\chi,x)=0=\hat{m}_{j}(\phi)$ for $1 \leq j \leq r-1,$ then the following holds
$$ w^r[(I_{w,p}^{\chi})(f,x)- f(x)]= \sum_{i=1}^{p} \frac{\beta_{i}}{i^{r}} \left( \frac{(\theta^{r}f)(x)}{r!} (m_{r}(\chi,x)+\hat{m}_{r}(\phi)) \right) + o (w^{-r}).$$
\end{Corollary}

\begin{Corollary}
Under the conditions of Theorem \ref{t4}, if moreover $\bar{M}_{k}^{p}(\chi) =0,\ $ for \\
$ k=1,2, \cdots (p-1),$ then for $f \in C^{(r)}(\mathbb{R}^+)$ with $r \geq p,$ we have
$$ \lim_{w \rightarrow \infty} w^{p}[(I_{w,p}^{\chi} f)(x) - f(x)] = \frac{(\theta^{p}f)(x)}{(p+1)!} \bar{M}_{p}^{p}(\chi,\phi).$$
\end{Corollary}
It is important to mention here that $M_{k}^{p}(\chi,\phi)$ does not vanish, in general. In order to have $M_{k}^{p}(\chi,\phi)=0$ for $k=0,1,2, \cdots p-1,$ we need to solve the following system
\begin{equation}
 \sum_{i=1}^{p} \beta_{i}=1,\  \sum_{i=1}^{p} \frac{\beta_{i}}{i}=0,\  \sum_{i=1}^{p} \frac{\beta_{i}}{i^{2}}=0,\ \cdots \sum_{i=1}^{p} \frac{\beta_{i}}{i^{p-1}}=0.
\end{equation}
The solution of the above system yields a linear combination which provides the convergence of order at least $p$ for functions $f \in C^{(p)}(\mathbb{R}^{+}).$
\par
\vspace{0.1cm}
In particular, let $f \in C^{(r)}(\mathbb{R}^{+}),\ r \geq 3.$ Then it is evident from Theorem 4.2 that the operator (\ref{main}) converges linearly, in general. In order to improve the rate of convergence, we need to solve the system (5.6) for $p=3$ to get the coefficients as $\beta_{1}=\frac{1}{2},\beta_{2}=-4$ and $\beta_{3}=\frac{9}{2}.$ Using these coefficient, we establish the following linear combination
\begin{equation} \label{combo2}
(I_{w,3}^{\chi,\phi})(f,x) = \frac{1}{2} (I_{w}^{\chi,\phi})(f,x)+ (-4) (I_{2w}^{\chi,\phi})(f,x) + \frac{9}{2} (I_{3w}^{\chi,\phi})(f,x)
\end{equation}
which corresponds to the asymptotic formula, for $r \geq 3,$ given by
$$ [(I_{w,3}^{\chi,\phi})(f,x) - f(x)]=  \sum_{j=3}^{r} \frac{(\theta^{j}f)(x)}{(j+1)! w^{j}} M_{j}^{3}(\chi,\phi) + o(w^{-r}).$$
It can be observe that the above linear combination produces the rate of convergence atleast 3 for $f \in C^{(r)}(\mathbb{R}^{+}),\ r \geq 3.$

\subsection{Quantitative estimates of linear combination}

\begin{thm}\label{t4}
Let $\chi,\phi$ be the kernels and $f \in \mathcal{C}^{(r)}(\mathbb{R}^+).$ Then we have
$$\Bigg| [(I_{w,p}^{\chi,\phi})(f,x)- f(x)] - \sum_{j=1}^{r} \frac{(\theta^{j} f)(x)}{j!\ w^{j}} \bar{M}_{j}^{p}(\chi,\phi) \Bigg| \leq \frac{2^r \|\theta^{r}(f-g)\|_{\infty}}{r! \ w^{r}} \sum_{i=1}^{p} \frac{\beta_{i}}{i^r} \hat{K} \left(\theta^{r}f,\frac{B}{A} \frac{1}{w(r+1)}   \right),$$
where $\displaystyle A=\sum_{i=1}^{p} \frac{\beta_{i}}{i^r} \left( M_{0}(\chi) \hat{M}_{r}(\phi)+ M_{r}(\chi) \hat{M}_{0}(\phi)\right)$ and $\displaystyle B=\sum_{i=1}^{p} \frac{\beta_{i}}{i^{r+1}} \left( M_{0}(\chi) \hat{M}_{r+1}(\phi)+  M_{r+1}(\chi) \hat{M}_{0}(\phi)\right).$
\end{thm}

\noindent\bf{Proof.}\rm \ In view of Theorem 5.1, we have \\

\noindent
$\displaystyle \Bigg| [(I_{w}^{\chi}f)(x) - f(x)] - \sum_{j=1}^{r} \frac{(\theta^{j} f)(x)}{j!\ w^{j}}\bar{M}_{j}^{p}(\chi,\phi) \Bigg|$
\begin{eqnarray*}
 = \Bigg| \sum_{i=1}^{p} \beta_{i} \sum_{k= - \infty}^{+\infty}\chi(e^{-k} x^{iw})\  (iw)  \int_{0}^{\infty} \phi(e^{-k} t^{iw}) h \left( \frac{t}{x} \right) (\log t- \log x)^r \ \frac{dt}{t}\Bigg|.\end{eqnarray*}

Now we substitute $\displaystyle R_{r}(f,x,t):= h \left( \frac{t}{x} \right) (\log t - \log x)^{r}.$ Then by using the estimate $$ \displaystyle |R_{r}(f,x,t)| \leq \frac{2}{r!} |\log t - \log x|^{r} \ \hat{K} \left( \theta^{r}f, \frac{|\log t- \log x|^{r}}{2(r+1)} \right),$$ we can write \\

\noindent $\displaystyle \Bigg| [(I_{w}^{\chi}f)(x) - f(x)] - \sum_{j=1}^{r} \frac{\theta^{j}f(x)}{j! \ w^{j}} \left( \sum_{\eta=0}^{j} \binom{j}{\eta} \hat{m}_{j-\eta}(\phi)\  m_{\eta}(\chi,x) \right) \Bigg|$
\begin{eqnarray*}
&\leq & \sum_{i=1}^{p} \beta_{i} \sum_{k= - \infty}^{+\infty} |\chi(e^{-k} x^{iw})| \ (iw) \int_{0}^{\infty}  |\phi(e^{-k} t^{iw})| \ \left(\frac{2}{r!} |\log t - \log x|^r \right) \ \hat{K} \left( \theta^{r}f, \frac{|\log t - \log x|}{2(r+1)} \right) \frac{dt}{t}\\
&\leq & \frac{2 \|\theta^{r}(f-g)\|_{\infty}}{r!} \sum_{i=1}^{p} \beta_{i} \sum_{k= - \infty}^{+\infty} |\chi(e^{-k} x^{iw})|  \ (iw) \int_{0}^{\infty}  |\phi(e^{-k} t^{iw})| \ |\log t - \log x|^r \frac{dt}{t} \ + \\&&
\frac{\|\theta^{r+1} g \|_{\infty}}{(r+1)!} \sum_{i=1}^{p} \beta_{i} \sum_{k= - \infty}^{+\infty} |\chi(e^{-k} x^{iw})|  \ (iw) \int_{0}^{\infty}  |\phi(e^{-k} t^{iw})| \ |\log t - \log x|^{r+1} \frac{dt}{t}\\
&:= & I_{1}+I_{2}.
\end{eqnarray*}
First we evaluate $I_{1}.$
\begin{eqnarray*}
|I_{1}| & \leq & \frac{2 \|\theta^{r}(f-g)\|_{\infty}}{r! \ w^{r}} \sum_{i=1}^{p} \frac{\beta_{i}}{i^r} \sum_{k= - \infty}^{+\infty} |\chi(e^{-k} x^{iw})| \int_{0}^{\infty}  |\phi(p)| |\log p+k-iw \log x|^r \ \frac{dp}{p} \\
&\leq & \frac{2^r \|\theta^{r}(f-g)\|_{\infty}}{r! w^{r}}\ \sum_{i=1}^{p} \frac{\beta_{i}}{i^r} \Bigg(\sum_{k= - \infty}^{+\infty} |\chi(e^{-k} x^{iw})| \int_{0}^{\infty}  |\phi(p)| |\log p|^r \ \frac{dp}{p}+ \\&&
\sum_{k= - \infty}^{+\infty} |\chi(e^{-k} x^{iw})| |k - iw \log x|^r \int_{0}^{\infty}  |\phi(p)| \ \frac{dp}{p} \Bigg)\\
& \leq & \frac{2^r \|\theta^{r}(f-g)\|_{\infty}}{r! \ w^{r}} \sum_{i=1}^{p} \frac{\beta_{i}}{i^r} \left( M_{0}(\chi) \hat{M}_{r}(\phi)+ M_{r}(\chi) \hat{M}_{0}(\phi) \right).
\end{eqnarray*}
Similarly, we obtain
\begin{eqnarray*}
|I_{2}| & \leq & \frac{\|\theta^{r+1} g \|_{\infty}}{(r+1)!\  w^{r+1}} \sum_{i=1}^{p} \frac{\beta_{i}}{i^{r+1}} \sum_{k= - \infty}^{+\infty} |\chi(e^{-k} x^{iw})|  \int_{0}^{\infty}  |\phi(p)|\ |\log p +k-iw \log x|^{r+1} \frac{dp}{p}\\
& \leq & \frac{2^{r}\|\theta^{r+1} g \|_{\infty}}{(r+1)!\ w^{r+1}} \sum_{i=1}^{p} \frac{\beta_{i}}{i^{r+1}} \left( M_{0}(\chi) \hat{M}_{r+1}(\phi)+  M_{r+1}(\chi) \hat{M}_{0}(\phi)\right).
\end{eqnarray*}
Now using the estimates of $I_{1}-I_{2}$ and then passing the infimum over $g \in \mathcal{C}^{(r+1)}(\mathbb{R}^+),$ we get the desired result.

\begin{Corollary}
Let $\chi,\phi$ be the kernel functions and $f \in \mathcal{C}^{(1)}(\mathbb{R}^+).$ Then we have
$$\Bigg| [(I_{w,p}^{\chi,\phi})(f,x)- f(x)] -  \frac{(\theta f)(x)}{w} \bar{M}_{1}^{p}(\chi,\phi) \Bigg| \leq \frac{2 \|\theta^{r}(f-g)\|_{\infty}}{w} \sum_{i=1}^{p} \frac{\beta_{i}}{i} \hat{K} \left(\theta f,\frac{B}{A} \frac{1}{2 w} \right),$$
where $\displaystyle A=\sum_{i=1}^{p} \frac{\beta_{i}}{i} \left( M_{0}(\chi) \hat{M}_{1}(\phi)+ M_{1}(\chi) \hat{M}_{0}(\phi)\right)$ and $\displaystyle B=\sum_{i=1}^{p} \frac{\beta_{i}}{i^{2}} \left( M_{0}(\chi) \hat{M}_{2}(\phi)+  M_{2}(\chi) \hat{M}_{0}(\phi)\right).$
\end{Corollary}

\section{Examples and graphical representation}
In this section, we present few examples of the kernels satisfying the assumptions of the discussed theory along with the graphical representation. We begin with the well-known Mellin B-spline kernel, see \cite{bardaro7,own,lnr}.

\subsection{Mellin B-splines kernel}
The Mellin B-splines of order $n$ are defined as
$$\bar{B}_{n}(x):= \frac{1}{(n-1)!} \sum_{j=0}^{n} (-1)^{j} {n \choose j} \bigg( \frac{n}{2}+\log x-j \bigg)_{+}^{n+1}\ , \hspace{1cm} x \in \mathbb{R}^+.$$
The Mellin transformation of $\bar{B}_{n}$ is given by
\begin{eqnarray} \label{splinefourier}
\hat{M}[\bar{B}_{n}](c+it) = \bigg( \frac{\sin(\frac{t}{2})}{(\frac{t}{2})} \Bigg)^{n},  \ \ \hspace{0.5cm} t \neq 0.
\end{eqnarray}
Since $\bar{B}_{n}(x) $ is compactly supported for every $ n \in \mathbb{N},$ the condition (K2) is satisfied. For (K1), we use the following lemma derived in \cite{own} using the Mellin Poisson summation formula \cite{butzer4}, which has the following form
$$ (i)^{j} \sum_{k= - \infty}^{+\infty} \chi(e^{-k} u) ( k-\log u)^{j} = \sum_{k= - \infty}^{+\infty} \hat{M}^{(j)}[\chi](2k \pi i) \ u^{-2 k \pi i}, \ \ \ \ \ \ \mbox{for} \ k \in \mathbb{Z}.$$

\begin{lem}\label{l1} \cite{own}
The condition $\displaystyle \sum_{k= - \infty}^{+\infty}  \chi(e^{-k} x^{w})=1$ is equivalent to the condition
\begin{equation*}
\hat{M}[\chi] (2k \pi i) =
     \begin{cases}
     {1,} &\quad\text{if} \ \ \ \ {k=0} \\
     {0,} &\quad\text{if} \ \ \ \ {\text{otherwise}}\\
   \end{cases}
\end{equation*}
Moreover, the condition $ m_{j}(\chi,u)=0$ for $ j= 1,2,...n$ is equivalent to the condition \\ $ \hat{M^{(j)}}[\chi](2k \pi i)= 0$ for $ j= 1,2...,n$ and $ \forall \ k \in \mathbb{Z}.$
\end{lem}
In view of Lemma \ref{l1} and (\ref{splinefourier}), it can be seen that $\bar{B}_{n}(x)$ satisfies the condition (K1) for every $n \in \mathbb{N}$ and $x \in \mathbb{R}^+.$
\vspace{0.3cm}

Now we consider $\chi(x):= \bar{B}_{4}(x)$ and $\phi(x):= \hat{B}_{2}(x).$ Using the Mellin Poisson summation formula, we obtain the algebraic moments for $\chi(x)$ as
$$ m_{0}(\chi)=1,\ m_{1}(\chi)=0, \ m_{2}(\chi)=\frac{1}{3}, \ m_{3}(\chi)=0.$$ Moreover, the algebraic moments for $\phi(x)$ are given by
$$\hat{m}_{0}(\phi)=1,\ \hat{m}_{1}(\phi)=0, \ \hat{m}_{2}(\phi)=\frac{1}{6}, \ \hat{m}_{3}(\phi)=0.$$ In view of Theorem 4.2, the asymptotic formula for $f \in C^{(3)}(\mathbb{R}^+)$ is given  by
$$ \lim_{w \rightarrow \infty} w^2[I_{w}^{\chi,\phi}f)(x) - f(x]= \frac{1}{4} (\theta^{2}f)(x).$$
But, the linear combination (\ref{combo2}) leads to the following asymptotic formula
$$\lim_{w \rightarrow \infty} w^3[I_{w,3}^{\chi,\phi}f)(x) - f(x]= 0 ,$$ which ensures the order of convergence atleast 3 for $f \in C^{(3)}(\mathbb{R}^+).$
\vspace{0.3cm}

Next we compare the approximation by the family of operators (\ref{main}) and the linear combinations of these operators using graphs and error estimates. It is apparent from \textit{Figure 1,2} and \textit{Table 1,2} that the linear combination of the proposed operators yield the better approximation.

\begin{figure}[h]
\centering
{\includegraphics[width=0.8 \textwidth]{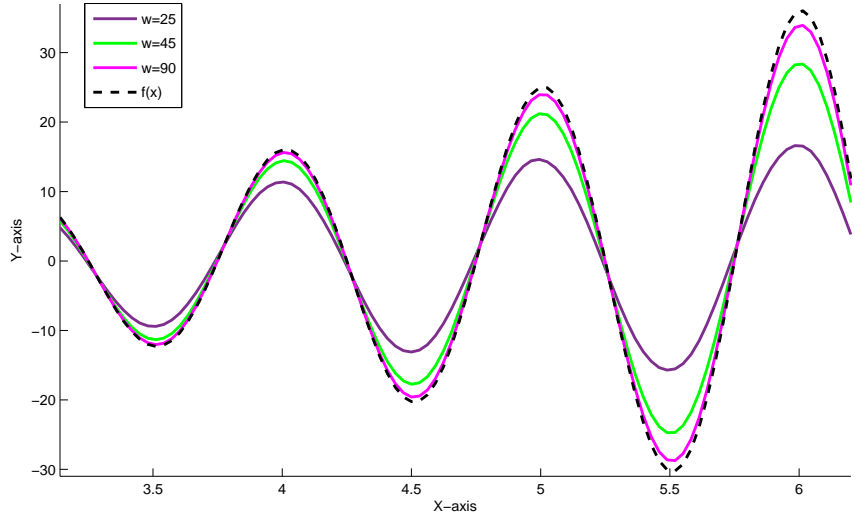}}
\caption{This figure shows the approximation of $ f(x)= x^{2} \cos(2 \pi x), \ x \in \left( \pi,2 \pi \right)$ by $(I_{w}^{\bar{B}_{4},\bar{B}_{4}}f)(x)$ for $w= 25, 45, 90$ respectively.}
\end{figure}

\begin{tb}\label{table1}\centering
 {\it Error estimation (upto $4$ decimal points) in the approximation of $f(x)$ by $(I_{w}^{\bar{B}_{4},\bar{B}_{4}}f)(x)$ for $w=25,45,90.$}

$  $

\begin{tabular}{|l|l|l|l|}\hline
 $x$&$ |f(x) - (I_{25}^{\bar{B}_{4},\bar{B}_{4}}f)(x)|$&$|f(x)-(I_{45}^{\bar{B}_{4},\bar{B}_{4}}f)(x)|$&$|f(x)-(I_{90}^{\bar{B}_{4},\bar{B}_{4}}f)(x)|$ \\
 \hline
 $3.55$&$2.9795$&$1.0033$&$0.2587$\\
  \hline
 $3.98$ & $4.4314$  &  $1.5008$ & $0.3876$ \\
  \hline
 $4.22$ & $1.7758$  &  $0.6922$ & $0.1869$\\
 \hline
  $4.85$ & $4.7175$  &  $1.5763$ & $0.4038$  \\
  \hline
 $5.35$&$6.7779$ & $2.3721$ & $0.6125$\\
 \hline

             \end{tabular}
   \end{tb}

\begin{figure}[h]
\centering
{\includegraphics[width=0.8 \textwidth]{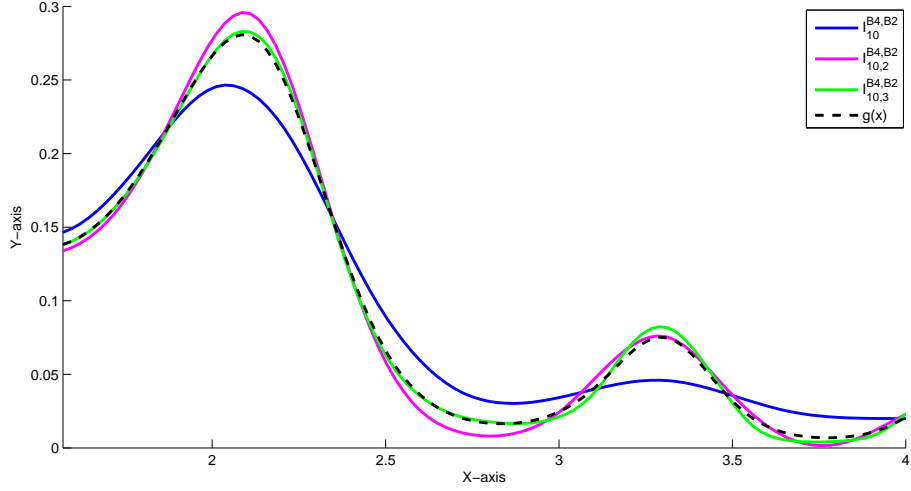}}
\caption{This figure exhibits the approximation of $g(x)= \frac{1}{x^3} e^{-\sin x^2},\ x \in \left( \frac{\pi}{2},4 \right) $ by $(I_{w}^{\bar{B}_{4},\bar{B}_{2}}g)(x)$ and  $(I_{w,i}^{\bar{B}_{4},\bar{B}_{2}}g)(x)$ for $i=2,3$  and $w=10$.}
\end{figure}

\begin{tb}\label{table2}\centering
 {\it Error estimation (upto $4$ decimal points) in the approximation of $g(x) $ by $(I_{w}^{\bar{B}_{4},\bar{B}_{2}}g)(x)$ and $(I_{w,i}^{\bar{B}_{4},\bar{B}_{2}}g)(x)$ for $i=2,3$  and $w=10.$}

$  $

\begin{tabular}{|l|l|l|l|}\hline
 $x$&$ |g(x) - (I_{10}^{\bar{B}_{4},\bar{B}_{2}}g)(x)|$&$|g(x)-(I_{10,2}^{\bar{B}_{4},\bar{B}_{2}}g)(x)|$&$|g(x)-(I_{10,3}^{\bar{B}_{4},\bar{B}_{2}}g)(x)|$ \\
 \hline
  $1.75$&$0.0087$&$0.0026$&$0.0007$\\
  \hline
 $2.10$&$0.0377$&$0.0153$&$0.0026$\\
  \hline
 $2.85$ & $0.0138$  &  $0.0076$ & $0.0002$ \\
  \hline
 $3.45$ & $0.0059$  &  $0.0037$ & $0.0021$\\
  \hline
 $3.95$&$0.0054$ & $0.0022$ & $0.0007$\\
 \hline
             \end{tabular}
   \end{tb}

\subsection{Combination of Mellin translates}
Consider the linear combination of translates of B-spline functions of order $n$ in the Mellin setting as follows:
\begin{eqnarray*}
\psi(x) &:=& c_{1} [\tau_{a} \bar{B}_{n}(x)] + c_{2}[ \tau_{b} \bar{B}_{n}(x)] \\
&=& c_{1}[\bar{B}_{n}(a x)]+  c_{2}[\bar{B}_{n}(b x)] ,\hspace{0.3cm \forall x \in \mathbb{R}^+,\ a,b \ \in \mathbb{R}.}
\end{eqnarray*}
Using the linearity of Mellin-tranform, we obtain
\begin{eqnarray} \nonumber \label{trans}
\hat{M}[\psi](w)&=& c_{1} \hat{M}[\bar{B}_{n}(a x)]+ c_{2} \hat{M} [\bar{B}_{n}(b x)]\\
&=& c_{1} \ a^{-w}  \hat{M} [\bar{B}_{n}](w)+ c_{2} \  b^{-w} \hat{M} [\bar{B}_{n}](w).
\end{eqnarray}
On differentiating (\ref{trans}), we have
$$ \hat{M}^{'}[\psi](w)= c_{1} (a^{-w} (\hat{M}^{'}[\bar{B}_{n}](w) - a^{-w} \log a) )+ c_{2} (b^{-w} (\hat{M}^{'}[\bar{B}_{n}](w) - b^{-w} \log b) ).$$
In view of Lemma \ref{l1}, we have the following system
$$ c_{1} + c_{2}=1 \ , \ \ c_{1} \log a + c_{2} \log b =0.$$
On solving for $c_{1} \ \ \mbox{and} \ \  c_{2},$ we get
$$ c_{1}= \frac{\log b}{(\log b - \log a)},\ \ c_{2}= \frac{- \log a}{ (\log b - \log a)}.$$
In particular, for $a= e^{-2} $ and $ b= e^{-3},$ we obtain the following linear combination of the Mellin's B-spline of order 2
$$ \psi(x) := 3 \bar{B}_{2}(e^{-2} x) - 2 \bar{B}_{2}(e^{-3} x)\ ,\ \ x \in \mathbb{R}^+.$$
Using Mellin Poisson summation formula, we obtain
$$m_{0}(\psi)= \hat{M}[\psi](0)= 1 \ \ \mbox{and} \ \ m_{1}(\psi)= \hat{M}^{'} [\psi](0)= 0.$$ Similarly we calculate
$m_{2}(\psi)= \frac{17}{3},\ m_{3}(\psi)= -32.$

\vspace{0.2cm}

Now for $f \in C^{(3)}(\mathbb{R}^+),$ the asymptotic formula is given as
$$ \lim_{w \rightarrow \infty} w^2[I_{w}^{\psi,\hat{B}_{2}}f)(x) - f(x]= \frac{35}{12} (\theta^{2}f)(x).$$
But, the linear combination (\ref{combo2}) leads to the following asymptotic formula:
$$\lim_{w \rightarrow \infty} w^3[I_{w}^{\psi,\hat{B}_{2}}f)(x) - f(x]= \frac{-16}{3} (\theta^{3}f)(x).$$ This shows that the combination (\ref{combo2}) provides order of convergence atleast 3 for $f \in C^{(3)}(\mathbb{R}^+).$

\end{document}